\documentclass{kybernetika-arXiv}
\usepackage{ucs}
\usepackage[utf8x]{inputenc}
\usepackage{amsmath}
\usepackage{amsfonts}
\usepackage{bbm}
\usepackage{bm}
\usepackage{amssymb}
\usepackage{amsthm}
\usepackage{mathtools}
\usepackage{mathrsfs}

\usepackage{natbib}
\newtheorem{theorem}{Theorem}[section]

\newtheorem{corollary}[theorem]{Corollary}

\newcommand{\eps}{\varepsilon}

\newcommand{\bY}{\mathbf{Y}}
\newcommand{\bZ}{\mathbf{Z}}
\newcommand{\bX}{\mathbf{X}}

\newcommand{\bI}{\mathbf{I}}

\newcommand{\bv}{\mathbf{v}}

\newcommand{\bt}{\mathbf{t}}

\newcommand{\zero}{\mathbf{0}}

\newcommand{\bbeta}{\bm{\beta}}
\newcommand{\beps}{\bm{\varepsilon}}

\newcommand{\btheta}{\bm{\Theta}}

\newcommand{\bdelta}{\bm{\Delta}}

\newcommand{\E}{\mathbbm{E}}

\newcommand{\Var}{\mathbb{V}\mbox{ar}\,}
\newcommand{\prob}{\mathbbm{P}}
\newcommand{\dist}{\mathscr{D}}

\newcommand{\Oo}{\mathcal{O}}

\begin{document}
\pagestyle{myheadings}

\title{Asymptotics for weakly dependent\\errors-in-variables}

% all authors
\author{Michal Pe\v{s}ta}

% contact informations of authors. Use one \contact comman per person.
% Specify name(s), surname(s), postal address and email.
\contact{Michal}{Pe\v{s}ta}{Charles University in Prague, Faculty of Mathematics and Physics, Department of Probability and Mathematical Statistics, Sokolovsk\'{a} 83, 186~75 Prague~8, Czech Republic}{pesta@karlin.mff.cuni.cz}

% heading contains authors and title of the paper, please shorten the title if necessary
% (it will appear only at the heading of every even page)
\markboth{Michal Pe\v{s}ta}{Asymptotics for weakly dependent errors-in-variables}

\maketitle

\begin{abstract}
Linear relations, containing measurement errors in input and output data, are taken into account in this paper. Parameters of these so-called \emph{errors-in-variables} (EIV) models can be estimated by minimizing the \emph{total least squares} (TLS) of the input-output disturbances. Such an~estimate is highly non-linear. Moreover in some realistic situations, the errors cannot be considered as independent by nature. \emph{Weakly dependent} ($\alpha$- and $\varphi$-mixing) disturbances, which are not necessarily stationary nor identically distributed, are considered in the EIV model. Asymptotic normality of the TLS estimate is proved under some reasonable stochastic assumptions on the errors. Derived asymptotic properties provide necessary basis for the validity of block-bootstrap procedures.
\end{abstract}

\keywords{errors-in-variables (EIV), dependent errors, total least squares (TLS), asymptotic normality}

\classification{15A51, 15A52, 62E20, 65F15, 62J99}

\section{INTRODUCTION}
The main goal is to establish results concerning asymptotic normality in linear relations, where weakly dependent measurement errors or disturbances in input and output data (errors-in-variables model) occur simultaneously.

\subsection{Errors-in-variables model}
\emph{Errors-in-variables} (EIV) model
\begin{equation}\label{EIV}\tag{E}
\mathop{\bY}_{n\times 1}=\mathop{\bZ}_{n\times p}\mathop{\bbeta}_{p\times 1}+\mathop{\beps}_{n\times 1}\quad\mbox{and}\quad\mathop{\bX}_{n\times p}=\mathop{\bZ}_{n\times p}+\mathop{\btheta}_{n\times p}
\end{equation}
is considered, where $\bbeta$ is a~vector of \emph{regression parameters} to be estimated, $\bX$ and $\bY$ consist of \emph{observable random} variables ($\bX$ are covariates and $\bY$ is a~response), $\bZ$ consists of \emph{unknown constants} and has full rank, and $\beps$ and $\btheta$ are \emph{random errors}.

\subsection{Weak Dependence}
In this paper, the EIV model is not restricted to independent observations due to the fact that in some situations the disturbances cannot be considered as independent by nature.% Hence, a~proper error structure is required and, consequently, suitable statistical inference needs to be derived.

%We are not in the case of independent observations any more and,
%Therefore, the dependence between measurement errors needs to be properly specified.
It is assumed that $\{\xi_n\}_{n=1}^{\infty}$ is a~sequence of random elements on a~probability space $(\Omega,\mathcal{F},\prob)$. For sub-$\sigma$-fields $\mathcal{A},\mathcal{B}\subseteq\mathcal{F}$, we define
\begin{align*}
\alpha(\mathcal{A},\mathcal{B})&:=\sup_{A\in\mathcal{A},B\in\mathcal{B}}\left|\prob(A\cap B)-\prob(A)\prob(B)\right|,\\
\varphi(\mathcal{A},\mathcal{B})&:=\sup_{A\in\mathcal{A},B\in\mathcal{B},\prob(A)>0}\left|\prob(B|A)-\prob(B)\right|.
\end{align*}
Intuitively, $\alpha$ and $\varphi$ measure the dependence of the events in $\mathcal{B}$ on those in $\mathcal{A}$. %Henceforth, let us define a~filtration $\mathcal{F}_m^n:=\sigma(\xi_i\in\mathcal{F},m\leq i\leq n)$.
There are many ways how to describe weak dependence or, in other words, \emph{asymptotic independence} of random variables \citep{Bradley2005}. We concentrate on two approaches. Considering a~filtration $\mathcal{F}_m^n:=\sigma(\xi_i\in\mathcal{F},m\leq i\leq n)$, sequence $\{\xi_n\}_{n=1}^{\infty}$ of random elements (e.g., variables) is said to be \emph{strong mixing}\index{strong mixing}\index{$\alpha$-mixing} ($\alpha$-mixing) if
\begin{equation}\label{alphamix}
\alpha(n):=\sup_{k\in\mathbbm{N}}\alpha(\mathcal{F}_{1}^k,\mathcal{F}_{k+n}^{\infty})\to 0,\quad n\to\infty;
\end{equation}
moreover, it is said to be \emph{uniformly strong mixing}\index{uniformly strong mixing}\index{$\varphi$-mixing} ($\varphi$-mixing) if
\begin{equation}\label{phimix}
\varphi(n):=\sup_{k\in\mathbbm{N}}\varphi(\mathcal{F}_{1}^k,\mathcal{F}_{k+n}^{\infty})\to 0,\quad n\to\infty.
\end{equation}
Uniformly strong mixing---introduced by~Rosenblatt---implies strong mixing~\citep{linlu1997}, which was proposed by~Ibragimov. Coefficients of dependence $\alpha(n)$ and $\varphi(n)$ measure how much dependence exists between events separated by at least $n$ observations or time periods.

\cite{Anderson1958} comprehensively analyzed a~class of $m$-dependent processes. They are $\varphi$-mixing, if they are finite order ARMA processes with innovations satisfying \emph{Doeblin's condition}\index{Doeblin's condition} \citep[p.~168]{Billingsley1968}. %, \citet[p.~192]{Doob1953}).
Finite order processes, which do not satisfy Doeblin's condition, can be shown to be $\alpha$-mixing \citep[pp.~312--313]{IL1971}. \cite{Rosenblatt1971} provides general conditions under which stationary Markov processes are $\alpha$-mixing. Since functions of mixing processes are themselves mixing \citep{Bradley2005}, time-varying functions of any of the processes just mentioned are mixing as well.

It is obvious that $\alpha(\mathcal{A},\mathcal{B})=\alpha(\mathcal{B},\mathcal{A})$ for arbitrary sub-$\sigma$-fields $\mathcal{A},\mathcal{B}\subseteq\mathcal{F}$. This symmetry does not hold for the $\varphi$-dependence. Indeed, \citet[pp.~213--214]{Rosenblatt1971} constructed some strictly stationary Markov chains that are $\varphi$-mixing but not ``time-reversed'' $\varphi$-mixing. Therefore, it is not possible to ``interchange'' the past with the future regarding the definition of the $\varphi$-mixing coefficient.

%The following lemma describes an~asymptotic behavior of $\alpha$- and $\varphi$-mixing coefficients of the corresponding random sequences after a~transformation. More precisely, the Borel transformation preserves the property of $\alpha$- and $\varphi$-mixing and, moreover, sustains the rate of the mixing coefficients.
%\begin{lemma}\label{lemma:bradley}
%Suppose that for each $m=1,2,\ldots$, $\xi^{(m)}:=\{\xi_k^{(m)}\}_{k\in\mathbbm{Z}}$ is a~sequence of random variables. Suppose the sequences $\xi^{(m)}$, $m=1,2,\ldots$ are independent of each other. Suppose that for each $k\in\mathbbm{Z}$, $h_k:\,\mathbbm{R}\times\mathbbm{R}\times\ldots\to\mathbbm{R}$ is a~Borel function. Define the sequence $\xi:=\{\xi_k\}_{k\in\mathbbm{Z}}$ of random variables by
%\[
%\xi_k:=h_k\left(\xi_k^{(1)},\xi_k^{(2)},\ldots\right),\quad k\in\mathbbm{Z}.
%\]
%Then for each $n\geq 1$, the following statements hold:
%\begin{enumerate}
%\item $\alpha(\xi,n)\leq\sum_{m=1}^{\infty}\alpha(\xi^{(m)},n)$,
%\item $\varphi(\xi,n)\leq\sum_{m=1}^{\infty}\varphi(\xi^{(m)},n)$.
%\end{enumerate}
%\end{lemma}

%\begin{Proof}
%See \citet[Theorem~5.2]{Bradley2005}.
%\end{Proof}

\subsection{Error Structure of the EIV model}\label{error-structure}
We think of all the vectors as columns. A~``row-column'' notation for a~matrix ${\bf A}$ is used in this manner: ${\bf A}_{i,\bullet}$ denotes the $i$-th row of matrix $\bf A$ and ${\bf A}_{\bullet,j}$ corresponds to the $j$-th column of matrix $\bf A$.

Let us consider a~probability space $(\Omega,\mathcal{F},\prob)$, where all the further mentioned random elements exist in. Proper distributional assumptions of random errors in the EIV model need to be proposed. Two levels of the error structure have to be distinguished. The first level of error structure---\emph{within-individual level}---is that each row $[\btheta_{i,\bullet},\eps_i]$ is absolutely continuous with respect to the Lebesgue measure having \emph{zero mean} and non-singular covariance matrix $\sigma^2\bI$, where $\sigma^2>0$ is unknown (for simplicity). This assumption can be straightforwardly generalized as discussed in, e.g., \cite{Pesta09robust} or \cite{Pesta09}. Relationships between individual observations are represented by the second level of error structure---\emph{between-individual level}. Here, the rows $[\btheta_{i,\bullet},\eps_i]$ form \emph{weakly dependent} sequences with \emph{zero mean}, i.e., zero mean $\alpha$-~or $\varphi$-mixing. The reason for this can come from the fact that the measurements, which are ``close to each other'', influence themselves somehow. Moreover, the influence decreases as the distance between observations increases.

Concerning the \emph{within-individual level}\index{error structure!within-individual level}, the mixing sequences of errors are assumed to be pairwise independent. The necessity and possible weakening of this assumption was discussed by~\citep{Pesta09robust}.

It has to be emphasized that any form of errors' \emph{stationarity} is not needed to assume. Omitting this, sometimes restrictive, assumption strengthen our results.

Additional \emph{design assumption} is necessary for asymptotics even in the case of independent errors:
\begin{equation}\tag{D}\label{delta}
\bdelta:=\lim_{n\to\infty}n^{-1}\bZ^{\top}\bZ\quad\mbox{exists and is positive definite}.
\end{equation}
Importance of the previous design assumption has already been thoroughly discussed in~\cite{Pesta09}.

\subsection{Total Least Squares Estimate}
Sometimes, full-information approaches like \emph{maximum likelihood} (ML) can provide parameter estimates for the previously mentioned model~\citep{healy}. Nevertheless, it is requisite to elaborate \emph{distributional-free estimation} method as well, e.g., total least squares introduced by~\cite{golub}.%, due to the impossibility to satisfy particular distributional assumptions in the model.

\emph{TLS estimate} $\widehat{\bbeta}$ of an~unknown parameter $\bbeta$ is defined as any vector $\bbeta$ satisfying~\eqref{EIV} such that the \emph{Frobenius norm} of error matrix $\left\|[\btheta,\beps]\right\|$ is minimal. Geometrically speaking, the Frobenius norm tries to minimize the \emph{orthogonal} distance between the observations and fitted hyperplane. For a~matrix ${\bf A}\equiv\left(A_{ij}\right)_{i,j=1}^{n,m}$, it is defined as
\begin{equation*}%\label{frobenius}
\left\|{\bf A}\right\|\equiv\left\|{\bf A}\right\|_F:=\sqrt{\sum_{i=1}^n\sum_{j=1}^m A_{ij}^2}=\sqrt{tr({\bf A}^{\top}{\bf A})}=\sqrt{\sum_{i=1}^{\min\{n,m\}} \sigma_i^2}=\sqrt{\sum_{i=1}^r \sigma_i^2},
\end{equation*}
where $r$ is the rank of matrix $\bf A$ and $\sigma_i$s are its singular values. The Frobenius norm can be viewed as a~\emph{multivariate} version of the Euclidean vector norm for matrices.

Let $\lambda$ be the $(p+1)$-st largest eigenvalue of a~matrix $[\bX,\bY]^{\top}[\bX,\bY]$ and let $\bv_{p+1}$ be the associated eigenvector. \citet{golub} showed that if $v_{p+1,p+1}\neq 0$, then the TLS estimate has form
\begin{equation*}%\label{TLSgallo}
\widehat{\bbeta}=(\bX^{\top}\bX-\lambda\bI)^{-1}\bX^{\top}\bY.
\end{equation*}
\citet{gleser} proved that with \emph{probability tending to one}, as $n$ increases, $v_{p+1,p+1}\neq 0$ and, hence, $\widehat{\bbeta}$ exists. In the sequel, it has to be kept in mind that $\lambda$ depends on $n$.

Finally, remind that the ML estimate of $\bbeta$ \citep{healy} \emph{coincides} with the TLS estimate if the rows of the error matrix are \emph{i.i.d.~multivariate normal}~with zero mean and non-singular covariance matrix.

\subsection{Consistency of the TLS estimate}
\emph{Strong consistency} of the TLS estimate for independent errors is proved by~\cite{gleser} and, moreover, \emph{weak consistency}---again for independent errors, but with less restrictive assumptions---is widely discussed in~\cite{gallo}. When a~premise of independence cannot be assumed, strong consistency of the TLS estimate under weak dependence of errors was explored by~\citep{Pesta09robust}.

Similar situation occurs to the TLS estimate's \emph{asymptotic normality}, which was proved by~\cite{gallophd} for the case of independent errors. This result is going to be extended for the weakly dependent errors now.

\section{Asymptotic Normality}
Firstly, central limit theorems for weakly dependent variables, that are going to be used in derivation of the TLS estimate's asymptotic normality, are needed to be stated. These \emph{CLTs cannot assume stationarity}, because they will not be applied directly on the (possibly stationary) errors of the EIV model, but on their transformations, which cannot be generally considered as stationary ones.
%A~\emph{weak invariance principle} (WIP) (also known as a~\emph{functional central limit theorem}) is a~functional convergence of the sum of variables to the standard Wiener process~$\mathcal{W}$. This principle for $\alpha$-mixing variables will be recalled.
Let us define $S_n:=\sum_{i=1}^n\xi_i$ and $\varsigma_n^2:=\Var S_n$. %Define random elements on Skorokhod space $D[0,1]$ as follows:
%\begin{equation*}
%\mathcal{W}_n(t):=\frac{S_{[nt]}}{\varsigma_n},\quad 0\leq t\leq 1,
%\end{equation*}
%where $[\cdot]$ denotes the nearest integer function. The expression $\varsigma_n^2$ is usually called \emph{long-run variance}\index{long-run variance}.
All proofs of the following results are in Appendix.

%\begin{lemma}[Weak invariance principle for $\alpha$-mixing]\label{lemma:WIPwALPHA}\index{weak invariance principle!$\alpha$-mixing}\index{WIP!$\alpha$-mixing}\index{central limit theorem!$\alpha$-mixing}
%Let $\{\xi_n\}_{n=1}^{\infty}$ be a~sequence of zero mean $\alpha$-mixing random variables with
%\begin{equation}\label{WIPwALPHAcond1}
%\sup_{n\in\mathbbm{N}}\E|\xi_n|^{2+\omega}<\infty
%\end{equation}
%and
%\begin{equation}\label{WIPwALPHAcond2}
%\sum_{n=1}^{\infty}\alpha(n)^{\omega/(2+\omega)}<\infty
%\end{equation}
%for some $\omega>0$. Suppose that
%\begin{equation}\label{WIPwALPHAcond3}
%\frac{\E S_n^2}{n}\to\varsigma^2>0,\quad n\to\infty
%\end{equation}
%is satisfied. Then
%\begin{equation*}%\label{WIPwALPHAresult}
%\mathcal{W}_n\xrightarrow{D[0,1]}\mathcal{W},\quad n\to\infty.\index{Wiener process}
%\end{equation*}
%\end{lemma}

%\begin{Proof}
%See \cite{Herrndorf1985} or \citet[Corollary~3.2.1]{linlu1997}.
%\end{Proof}

%Since the central limit theorem can be derived from the weak invariance principle, then a~corollary of previous Lemma~\ref{lemma:WIPwALPHA} can be stated.

\begin{corollary}[Central limit theorem for $\alpha$-mixing]\label{col:CLTwALPHA}\index{central limit theorem!$\alpha$-mixing}\index{CLT!$\alpha$-mixing}
%Suppose that all the assumptions of Lemma~\ref{lemma:WIPwALPHA} on a~sequence of zero mean $\alpha$-mixing random variables $\{\xi_n\}_{n=1}^{\infty}$ are satisfied.
Let $\{\xi_n\}_{n=1}^{\infty}$ be a~sequence of zero mean $\alpha$-mixing random variables with
\begin{equation}\label{WIPwALPHAcond1}
\sup_{n\in\mathbbm{N}}\E|\xi_n|^{2+\omega}<\infty
\end{equation}
and
\begin{equation}\label{WIPwALPHAcond2}
\sum_{n=1}^{\infty}\alpha(n)^{\omega/(2+\omega)}<\infty
\end{equation}
for some $\omega>0$. Suppose that
\begin{equation}\label{WIPwALPHAcond3}
\frac{\E S_n^2}{n}\to\varsigma^2>0,\quad n\to\infty
\end{equation}
is satisfied. Then
$\frac{S_n}{\varsigma_n}\xrightarrow{\dist}\mathcal{N}(0,1),\,n\to\infty.$
\end{corollary}

%Similar CLT for $\varphi$-mixing sequences is needed as well.

%\begin{lemma}[Lindeberg central limit theorem for $\varphi$-mixing]\label{lemma:CLTwPHIlin}\index{central limit theorem!$\varphi$-mixing!Lindeberg}\index{CLT!$\varphi$-mixing!Lindeberg}
%Let $\{\xi_n\}_{n=1}^{\infty}$ be a~sequence of zero mean $\varphi$-mixing random variables having finite variance. Suppose that the Lindeberg condition
%\begin{equation}\label{CLTwPHIlincond}
%\forall\delta>0:\quad\lim_{n\to\infty}\frac{1}{\varsigma_n^2}\sum_{i=1}^n\E\xi_i^2\mathcal{I}\{|\xi_i|>\delta\varsigma_n\}=0
%\end{equation}
%is satisfied. Then
%\begin{equation*}
%\frac{S_n}{\varsigma_n}\xrightarrow{\dist}\mathcal{N}(0,1),\quad n\to\infty.
%\end{equation*}
%\end{lemma}

%\begin{Proof}
%See \citet[Corollary~4]{Utev1990}.
%\end{Proof}

%Lindeberg condition~\eqref{CLTwPHIlincond} can be replaced by a~stronger type of Lyapunov condition. This fact leads into the following corollary, which is more comfortable for us from the point of applicability.

\begin{corollary}[Central limit theorem for $\varphi$-mixing]\label{col:CLTwPHI}
Let $\{\xi_n\}_{n=1}^{\infty}$ be a~sequence of zero mean $\varphi$-mixing random variables such that
\begin{equation}\label{CLTwPHIcond1}
\sup_{n\in\mathbbm{N}}\E|\xi_n|^{2+\omega}<\infty
\end{equation}
for some $\omega>0$ and
\begin{equation}\label{CLTwPHIcond2}
\frac{\E S_n^2}{n}\to\varsigma^2>0,\quad n\to\infty.
\end{equation}
Then $\frac{S_n}{\varsigma_n}\xrightarrow{\dist}\mathcal{N}(0,1),\, n\to\infty.$
\end{corollary}

%Assumption~\eqref{CLTwPHIcond2} may even be replaced by a~weaker one:
%\[
%\liminf_{n\to\infty}\frac{\E S_n^2}{n}=\varsigma^2>0,
%\]
%where the limit inferior is used instead of the original limit.

For a~given sequence $\xi_{\circ}\equiv\{\xi_n\}_{n=1}^{\infty}$ of random elements, the dependence coefficients $\alpha(n)$ is denoted $\alpha(\xi_{\circ},n)$. Analogous notation is used for $\varphi$-mixing sequences.

The first main result of this section is the~asymptotic normality for the TLS estimate, where the errors are $\alpha$-mixing.

\begin{theorem}[Asymptotic normality in $\alpha$-mixing EIV]\label{thm:ANalpha}
Let the EIV model hold and assumption~\eqref{delta} be satisfied. Suppose $\{\Theta_{n,1}\}_{n=1}^{\infty},\ldots,\{\Theta_{n,p}\}_{n=1}^{\infty}$, and $\{\eps_n\}_{n=1}^{\infty}$ are pairwise independent sequences of $\alpha$-mixing random variables having
\begin{equation}\label{assalpha1AN}
\alpha(\Theta_{\circ,j},n)=\Oo(n^{-1-\delta_j}),\, j=1,\ldots,p
\quad\mbox{and}\quad
\alpha(\eps_{\circ},n)=\Oo(n^{-1-\delta_{p+1}}),
\end{equation}
as $n\to\infty$ for some $\delta_j>0$, $j\in\{1,\ldots,p+1\}$. Moreover, assume that
\begin{equation}\label{assalphaZAN}
\sup_{n\in\mathbbm{N}} Z_{n,j}^2<\infty,\, j\in\{1,\ldots,p\},
\end{equation}
\begin{equation}\label{assalpha3AN1}
\sup_{n\in\mathbbm{N}}\E|\Theta_{n,j}|^{4+\omega_j}<\infty,\, j\in\{1,\ldots,p\},
\quad\mbox{and}\quad
\sup_{n\in\mathbbm{N}}\E|\eps_n|^{4+\omega_{p+1}}<\infty
\end{equation}
for some $\omega_j>0$, $j\in\{1,\ldots,p+1\}$ such that
\begin{equation}\label{assalpha4AN}
\frac{2}{\min_{j=1,\ldots,p+1}\omega_j} < \min_{j=1,\ldots,p+1}\delta_j.
\end{equation}
If there exists a~positive definite matrix $\beth$ such that
\begin{equation}\label{long-run}
n^{-1}\Var[\bX,\bY]^{\top}[\bX,\bY]\left[\begin{array}{c}
\bbeta\\-1
\end{array}\right]\to\beth>\zero,\quad n\to\infty;
\end{equation}
then $\sqrt{n}\left(\widehat{\bbeta}-\bbeta\right)\xrightarrow{\dist}\mathcal{N}(\zero,\cdot),\, n\to\infty.$
\end{theorem}

The covariance matrix of the normal distribution is indeed omitted. The reason is that even in the simpler case of independent errors, three main pitfalls exist~\citep{Pesta09}: the covariance of the limiting multivariate normal distribution depends on the unknown parameter~$\bbeta$ and on the unknown matrix~$\bdelta$, and without the assumption on the third and fourth moments of the rows of~$[\btheta,\beps]$, the covariance matrix has a~very complicated form. A~partial solution to the first two mentioned issues could be plugging consistent estimates instead of the unknown entities. On the contrary, the third issue seems to be a~big problem whatsoever, because the third and the fourth errors' moments cannot be estimated from the data.

Therefore, it is a~waste of effort to theoretically derive the covariance, because such a~formula---which could be unobtainable---is already computationally useless. It is sufficient to know that the TLS estimate has a~non-degenerate asymptotic normal distribution with some (unknown) positive definite covariance matrix. Nevertheless, the practical calculation of the asymptotic distribution will be proceeded by the resampling methods anyway.

The second main result of this section is the asymptotic normality of the TLS estimate, where the errors are $\varphi$-mixing.

\begin{theorem}[Asymptotic normality in $\varphi$-mixing EIV]\label{thm:ANphi}
Let the EIV model hold and assumption~\eqref{delta} be satisfied. Suppose $\{\Theta_{n,1}\}_{n=1}^{\infty},\ldots,\{\Theta_{n,p}\}_{n=1}^{\infty}$, and $\{\eps_n\}_{n=1}^{\infty}$ are pairwise independent sequences of $\varphi$-mixing random variables such that
\begin{equation}\label{assphi1ANa}
\sum_{n=1}^{\infty}\sqrt{\varphi(\Theta_{\circ,j},n)}<\infty,\, j\in\{1,\ldots,p\}
\quad\mbox{and}\quad
\sum_{n=1}^{\infty}\sqrt{\varphi(\eps_{\circ},n)}<\infty.
\end{equation}
Moreover, assume that
\begin{equation}\label{assphiZAN}
\sup_{n\in\mathbbm{N}} Z_{n,j}^2<\infty,\, j\in\{1,\ldots,p\},
\end{equation}
\begin{equation}\label{assphi2AN1}
\sup_{n\in\mathbbm{N}}\E|\Theta_{n,j}|^{4+\omega_j}<\infty,\, j\in\{1,\ldots,p\},
\quad\mbox{and}\quad
\sup_{n\in\mathbbm{N}}\E|\eps_n|^{4+\omega_{p+1}}<\infty
\end{equation}
for some $\omega_j>0$, $j\in\{1,\ldots,p+1\}$. If there exists a~positive definite matrix $\beth$ such that
\begin{equation}\label{long-run2}
n^{-1}\Var[\bX,\bY]^{\top}[\bX,\bY]\left[\begin{array}{c}
\bbeta\\-1
\end{array}\right]\to\beth>\zero,\quad n\to\infty;
\end{equation}
then $\sqrt{n}\left(\widehat{\bbeta}-\bbeta\right)\xrightarrow{\dist}\mathcal{N}(\zero,\cdot),\, n\to\infty.$
\end{theorem}

\section{Conclusions}
A~structure of the EIV model with \emph{weakly dependent errors} is introduced in this paper. Suitable central limit theorems for strong mixing and uniformly strong mixing sequences are postulated. \emph{Strong consistency} of the TLS estimate under both proposed forms of errors' \emph{asymptotic independence} is used for proving \emph{asymptotic normality} of the TLS estimate under $\alpha$- and $\varphi$-mixing errors in the EIV model. Stochastic assumptions for the asymptotic normality of the TLS estimate are derived. Moreover, all the presented results are valid for non-stationary errors.
 
%Similar situation as in the independent errors' case would occur with the calculation of variance of the TLS estimate. Since the variance is already computationally useless when the errors are independent, another approach, that provides the approximate (estimated) asymptotic variance of the TLS estimate, needs to be proposed.

\subsection{Discussion}
A~straightforward application of the derived asymptotic properties lies in the \emph{justification of the bootstrap methods}. The variance of the asymptotic multivariate normal distribution of the TLS estimate could sometimes be computationally unfeasible. One possible solution to this dilemma is bootstrapping. \cite{Pesta09} proved validity of various bootstrap procedures in the case of independent errors in the EIV model, where asymptotic normality of the TLS estimate served as a~crucial cornerstone. Asymptotic results from this paper will provide a~workable basis for the block-bootstrap techniques when constructing confidence intervals and testing hypotheses in the case of weakly dependent errors of the EIV model.

It is known that a~sequence of random variables is $\varphi$-mixing implies that this sequence is $\alpha$-mixing. On the other hand, the CLT for $\varphi$-mixing (Corollary~\ref{col:CLTwPHI}) has weaker assumptions than the CLT for $\alpha$-mixing (Corollary~\ref{col:CLTwALPHA}). Indeed, Corollary~\ref{col:CLTwPHI} does not require any assumptions on \emph{mixing rate} $\varphi(n)$ such as assumption~\eqref{WIPwALPHAcond2} on $\alpha$-mixing rates. Therefore in Theorem~\ref{thm:ANphi}, we do not have to deal with mixing rate assumption like~\eqref{assalpha1AN} nor a~restriction on the \emph{moment order} like~\eqref{assalpha4AN} as in Theorem~\ref{thm:ANalpha}. On the contrary, Theorem~\ref{col:CLTwPHI} requires mixing rate assumption~\eqref{assphi1ANa}, which are inherited from the assumptions for the strong consistency of the TLS estimate. Both asymptotic normality results of the TLS estimate (Theorem~\ref{thm:ANalpha} and Theorem~\ref{thm:ANphi}) require all the assumptions from the strong consistency results (Theorem~\ref{thm:CONalpha} and Theorem~\ref{thm:CONphi}), because they were used in the proofs of asymptotic normality.

Assumptions~\eqref{long-run} and~\eqref{long-run2} concerning the \emph{long-run variance} of the TLS estimate are requisite and cannot be omitted, because they assure that the variance of the TLS estimate is bounded away from zero and, simultaneously, does not explode into infinity. These assumptions straightforwardly allow to apply the appropriate CLT in order to prove the asymptotic normality of $\sqrt{n}(\widehat{\bbeta}-\bbeta)$.

All the assumptions and remarks regarding the error structure or the pairwise independence of errors have been already discussed by~\cite{Pesta09robust}.

Finally, if identically distributed rows of errors (on the between-individual level) are taken into account together with finiteness of their moments up to a~suitable order, then moments' assumptions~\eqref{assalpha3AN1} and~\eqref{assphi2AN1} are trivially satisfied. Then, a~\emph{strict stationarity} of the between-individual errors and existence of appropriate moments have to imply these assumptions as well. In other words, moment assumptions~\eqref{assalpha3AN1} and~\eqref{assphi2AN1} cannot be considered as unattainable. Moreover for strictly stationary errors, even the supremum in definitions~\eqref{alphamix} and~\eqref{phimix} can be simply avoided.

%Zero mean errors are implicitly assumed through the whole paper and are not explicitly specified in every theorem, because this assumption is a~part of the EIV model's definition valid for the whole paper.

\section*{ACKNOWLEDGEMENT}
\small
This paper was written with the support of the Czech Science Foundation project ``DYME – Dynamic Models in Economics'' No.~P402/12/G097.

\makesubmdate

\appendix

\section{Appendix}

\begin{theorem}[Strong consistency in $\alpha$-mixing EIV]\label{thm:CONalpha}
Let the EIV model hold and assumption~\eqref{delta} be satisfied. Suppose $\{\Theta_{n,1}\}_{n=1}^{\infty},\ldots,\{\Theta_{n,p}\}_{n=1}^{\infty}$, and $\{\eps_n\}_{n=1}^{\infty}$ are pairwise independent sequences of $\alpha$-mixing random variables having
\begin{equation}\label{assalpha1}
\alpha(\Theta_{\circ,j},n)=\Oo(n^{-q_j/(2q_j-2)-\delta_j}),\, j=1,\ldots,p
\quad\mbox{and}\quad
\alpha(\eps_{\circ},n)=\Oo(n^{-q_{p+1}/(2q_{p+1}-2)-\delta_{p+1}}),
\end{equation}
as $n\to\infty$ for some $\delta_j>0$ and $1<q_j\leq 2$, $j\in\{1,\ldots,p+1\}$. If
\begin{equation}\label{assalphaZ}
\sup_{n\in\mathbbm{N}} Z_{n,j}^2<\infty,\, j\in\{1,\ldots,p\},
\end{equation}
\begin{equation}\label{assalpha3a}
\sup_{n\in\mathbbm{N}}\E|\Theta_{n,j}|^{2q_j}<\infty,\, j\in\{1,\ldots,p\},
\quad\mbox{and}\quad
\sup_{n\in\mathbbm{N}}\E|\eps_n|^{2q_{p+1}}<\infty,
\end{equation}
then
\begin{equation}\label{alphaRESULT1}
\lim_{n\to\infty}\widehat{\bbeta}=\bbeta,\,\mbox{a.s.}
\quad\mbox{and}\quad
\lim_{n\to\infty}\frac{\lambda}{n}=\sigma^2,\,\mbox{a.s.}
\end{equation}
\end{theorem}

\begin{Proof}
See~\citet[Theorem~3.1]{Pesta09robust}.
\end{Proof}

\begin{theorem}[Strong consistency in $\varphi$-mixing EIV]\label{thm:CONphi}
Let the EIV model hold and assumption~\eqref{delta} be satisfied. Suppose $\{\Theta_{n,1}\}_{n=1}^{\infty},\ldots,\{\Theta_{n,p}\}_{n=1}^{\infty}$, and $\{\eps_n\}_{n=1}^{\infty}$ are pairwise independent sequences of $\varphi$-mixing random variables such that
\begin{equation}\label{assphi1a}
\sum_{n=1}^{\infty}\sqrt{\varphi(\Theta_{\circ,j},n)}<\infty,\, j\in\{1,\ldots,p\}
\quad\mbox{and}\quad
\sum_{n=1}^{\infty}\sqrt{\varphi(\eps_{\circ},n)}<\infty.
\end{equation}
If
\begin{equation}\label{assphi2a}
\sum_{n=1}^{\infty}\frac{\E\Theta_{n,j}^4}{n^2}<\infty,\, j\in\{1,\ldots,p\}
\quad\mbox{and}\quad
\sum_{n=1}^{\infty}\frac{\E\eps_n^4}{n^2}<\infty,
\end{equation}
then
\begin{equation}\label{phiRESULT1}
\lim_{n\to\infty}\widehat{\bbeta}=\bbeta,\,\mbox{a.s.}\quad\mbox{and}\quad
\lim_{n\to\infty}\frac{\lambda}{n}=\sigma^2,\,\mbox{a.s.}
\end{equation}
\end{theorem}

\begin{Proof}
See~\citet[Theorem~3.2]{Pesta09robust}.
\end{Proof}

\begin{Proof}[\emph{Corollary~\ref{col:CLTwALPHA}}]

Define random elements on Skorokhod space $D[0,1]$ as follows:
\begin{equation*}
\mathcal{W}_n(t):=\frac{S_{[nt]}}{\varsigma_n},\quad 0\leq t\leq 1,
\end{equation*}
where $[\cdot]$ denotes the nearest integer function. Then weak invariance principle for $\alpha$-mixing by~\cite{Herrndorf1985} or \citet[Corollary~3.2.1]{linlu1997} completes this proof, because a~functional distributional limit (a~convergence in Skorokhod space) on a~bounded interval implies the pointwise distributional limit at the right end-point of the interval. Indeed, there is a~problem with measurability of the functionals of discontinuous processes. Space of c\`{a}dl\`{a}g functions with the uniform metric is not separable, but equipped with the Skorokhod metric is~\citep[Section~14]{Billingsley1968}. Recall the Skorokhod metric $\varrho_{[0,1]}(\cdot,\cdot)$ for c\`{a}dl\`{a}g functions on $[0,1]$ by introducing $\mathcal{F}_{[0,1]}:=\left\{f:[0,1]\xleftrightarrow{1-1} [0,1],\,\mbox{strictly increasing}\right\}$ and, hence, $\varrho_{[0,1]}(g,h):=\inf_{f\in\mathcal{F}_{[0,1]}}\max\left\{\sup_{t\in[0,1]}|f(t)-1|,\sup_{t\in[0,1]}|g(t)-h(f(t))|\right\}$. If $f\in\mathcal{F}_{[0,1]}$, then $f(1)=1$. Since
\[
\max\left\{\sup_{t\in[0,1]}|f(t)-1|,\sup_{t\in[0,1]}|\mathcal{W}_n(t)-\mathcal{W}(f(t))|\right\}\geq |\mathcal{W}_n(1)-\mathcal{W}(1)|,
\]
where~$\mathcal{W}$ is the standard Wiener process, then
\begin{multline*}
0\xleftarrow{n\to\infty}\varrho_{[0,1]}(\mathcal{W}_n,\mathcal{W})=\inf_{f\in\mathcal{F}_{[0,1]}}\max\left\{\sup_{t\in[0,1]}|f(t)-1|,\sup_{t\in[0,1]}|\mathcal{W}_n(t)-\mathcal{W}(f(t))|\right\}\\\geq |\mathcal{W}_n(1)-\mathcal{W}(1)|.
\end{multline*}
\end{Proof}

\begin{Proof}[\emph{Corollary~\ref{col:CLTwPHI}}]

It is sufficient to show that assumptions~\eqref{CLTwPHIcond1} and~\eqref{CLTwPHIcond2} implies Lindeberg condition
\[
\forall\delta>0:\quad\lim_{n\to\infty}\frac{1}{\varsigma_n^2}\sum_{i=1}^n\E\xi_i^2\mathcal{I}\{|\xi_i|>\delta\varsigma_n\}=0
\]
from \citet[Corollary~4]{Utev1990}. The first step is to show that conditions~\eqref{CLTwPHIcond1} and~\eqref{CLTwPHIcond2} implies so-called Lyapunov condition, i.e., having fixed $\omega>0$:
\[
\frac{1}{\varsigma_n^{2+\omega}}\sum_{i=1}^n\E|\xi_i|^{2+\omega}\leq\frac{1}{\varsigma_n^{2+\omega}}\sum_{i=1}^n\sup_{\iota\in\mathbbm{N}}\E|\xi_{\iota}|^{2+\omega}=\frac{n}{\varsigma_n^{2+\omega}}\sup_{\iota\in\mathbbm{N}}\E|\xi_{\iota}|^{2+\omega}\to 0,\quad n\to\infty.
\]

Now, Lyapunov condition $\lim_{n\to\infty}\varsigma_n^{-2-\omega}\sum_{i=1}^n\E|\xi_i|^{2+\omega}=0$ holds and we fix $\delta>0$. Since $|\xi_i|>\delta\varsigma_n$ implies $|\xi_i/\delta\varsigma_n|^{\omega}>1$, we obtain
\begin{multline*}
\frac{1}{\varsigma_n^2}\sum_{i=1}^n\E\xi_i^2\mathcal{I}\{|\xi_i|>\delta\varsigma_n\}\leq\frac{1}{\delta^{\omega}\varsigma_n^{2+\omega}}\sum_{i=1}^n\E|\xi_i|^{2+\omega}\mathcal{I}\{|\xi_i|>\delta\varsigma_n\}\\
\leq\frac{1}{\delta^{\omega}\varsigma_n^{2+\omega}}\sum_{i=1}^n\E|\xi_i|^{2+\omega}\to 0,\quad n\to\infty.
\end{multline*}
\end{Proof}

\begin{Proof}[\emph{Theorem~\ref{thm:ANalpha}}]

Assumptions of Theorem~\ref{thm:ANalpha} imply the assumptions of Theorem~\ref{thm:CONalpha}. Indeed, assumptions~\eqref{assalphaZ} and~\eqref{assalphaZAN} coincide. Assumption~\eqref{assalpha1AN} on $\alpha$-mixing rates clearly implies assumption~\eqref{assalpha1} for any $\delta_j>0$ and $1<q_j\leq 2$, $j\in\{1,\ldots,p+1\}$. Supremum assumption~\eqref{assalpha3AN1} implies~\eqref{assalpha3a} for any $\omega_j>0$ and $1<q_j\leq 2$, $j\in\{1,\ldots,p+1\}$ as well, because of a~corollary of the Jensen's inequality
\begin{equation}\label{jensen-moment}
(\E|\xi|^r)^{1/r}\leq(\E|\xi|^s)^{1/s},\,0<r<s<\infty.
\end{equation}

\citet[Section~2]{golub} proved that $[\widehat{\boldsymbol\beta},-1]^{\top}$ is the eigenvector of matrix $\left[\bX,\bY\right]^{\top}\left[\bX,\bY\right]$ with associated eigenvalue $\lambda$. Hence,
\begin{equation*}%\label{closed-derivation}
\left[\bX,\bY\right]^{\top}\left[\bX,\bY\right]\left[\begin{array}{c}
\widehat{\boldsymbol\beta}\\-1
\end{array}\right]=\left[\begin{array}{cc}
{\bf X}^{\top}{\bf X}&{\bf X}^{\top}{\bf Y}\\
{\bf Y}^{\top}{\bf X}&{\bf Y}^{\top}{\bf Y}
\end{array}\right]\left[\begin{array}{c}
\widehat{\boldsymbol\beta}\\-1
\end{array}\right]=\lambda\left[\begin{array}{c}
\widehat{\boldsymbol\beta}\\-1
\end{array}\right].
\end{equation*}
Previous partitioning yields $\bX^{\top}\bX\widehat{\bbeta}-\bX^{\top}\bY=\lambda\widehat{\bbeta}$ and $\bY^{\top}\bX\widehat{\bbeta}-\bY^{\top}\bY=-\lambda$, which can be rewritten in the following manner $\bX^{\top}\bY=(\bX^{\top}\bX-\lambda\bI)\widehat{\bbeta}$ and $\bY^{\top}\bY=\bY^{\top}\bX\widehat{\bbeta}+\lambda$. Hence,
\begin{align*}
[\bX,\bY]^{\top}[\bX,\bY]&=\left[\begin{array}{cc}
\bX^{\top}\bX & (\bX^{\top}\bX-\lambda\bI)\widehat{\bbeta}\\
\widehat{\bbeta}^{\top}(\bX^{\top}\bX-\lambda\bI) & \bY^{\top}\bX\widehat{\bbeta}+\lambda
\end{array}\right]\\&=\left[\bI,\widehat{\bbeta}\right]^{\top}(\bX^{\top}\bX-\lambda\bI)\left[\bI,\widehat{\bbeta}\right]+\lambda\bI.
\end{align*}
It can be easily noted that $\left[\bI,\bbeta\right]\left[\begin{array}{c}
\bbeta\\-1
\end{array}\right]=\zero$. Therefore, we obtain
\begin{equation*}%\label{starting}
\left[\bI,\bbeta\right][\bX,\bY]^{\top}[\bX,\bY]\left[\begin{array}{c}
\bbeta\\-1
\end{array}\right]=\left[\bI,\bbeta\right]\left[\bI,\widehat{\bbeta}\right]^{\top}(\bX^{\top}\bX-\lambda\bI)(\bbeta-\widehat{\bbeta})
\end{equation*}
and, then,
\begin{equation}\label{hatbetaminus}
\sqrt{n}(\widehat{\bbeta}-\bbeta)=-\bdelta_n^{-1}\left(\left[\bI,\bbeta\right]\left[\bI,\widehat{\bbeta}\right]^{\top}\right)^{-1}\left[\bI,\bbeta\right]\left(n^{-1/2}[\bX,\bY]^{\top}[\bX,\bY]\right)\left[\begin{array}{c}
\bbeta\\-1
\end{array}\right].
\end{equation}
Since $\E[\bX,\bY]^{\top}[\bX,\bY]=[\bI,\bbeta]^{\top}\bZ^{\top}\bZ[\bI,\bbeta]+n\sigma^2\bI$, then
\begin{equation}\label{Ezero}
[\bI,\bbeta]\E[\bX,\bY]^{\top}[\bX,\bY]\left[\begin{array}{c}
\bbeta \\ -1
\end{array}\right]=\zero.
\end{equation}
Relation~\eqref{hatbetaminus} can be alternatively rewritten using identity~\eqref{Ezero} in a~slightly more sophisticated way
\begin{multline*}%\label{betahattoprove}
\sqrt{n}(\widehat{\bbeta}-\bbeta)=-\bdelta_n^{-1}\left(\left[\bI,\bbeta\right]\left[\bI,\widehat{\bbeta}\right]^{\top}\right)^{-1}\left[\bI,\bbeta\right]\\\left(n^{-1/2}\left\{[\bX,\bY]^{\top}[\bX,\bY]-\E[\bX,\bY]^{\top}[\bX,\bY]\right\}\right)\left[\begin{array}{c}
\bbeta\\-1
\end{array}\right],
\end{multline*}
which will be useful in the forthcoming steps of this proof.

The proof of Theorem~3.1 by~\cite{Pesta09robust} also provides the following important relation
\begin{equation}\label{gle1}
\bdelta_n:=n^{-1}(\bX^{\top}\bX-\lambda\bI)\xrightarrow[n\to\infty]{a.s.}\bdelta.
\end{equation}

%Assumptions of this theorem imply the assumptions of Theorem~\ref{thm:CONalpha} and, hence, its consistency results can be used.
With respect to~\eqref{gle1} and Theorem~\ref{thm:CONalpha}, we have
\begin{equation}\label{betaslucky}
\bdelta_n^{-1}\left(\left[\bI,\bbeta\right]\left[\bI,\widehat{\bbeta}\right]^{\top}\right)^{-1}\left[\bI,\bbeta\right]\xrightarrow{a.s.}\bdelta^{-1}\left(\left[\bI,\bbeta\right]\left[\bI,\bbeta\right]^{\top}\right)^{-1}\left[\bI,\bbeta\right],\quad n\to\infty.
\end{equation}
With probability tending to one, the inverse $\bdelta_n$ exists due to~\eqref{delta} and~\eqref{gle1}, and the inverse of $\left[\bI,\bbeta\right]\left[\bI,\widehat{\bbeta}\right]^{\top}$ exists because of~\eqref{alphaRESULT1}. Moreover, matrix $\left[\bI,\bbeta\right]\left[\bI,\bbeta\right]^{\top}=\bI+\bbeta\bbeta^{\top}$ is always positive definite and, hence, regular.

Convergence almost surely from~\eqref{betaslucky} and the Slutsky's theorem reduce the problem of finding a~limiting distribution for $\sqrt{n}(\widehat{\bbeta}-\bbeta)$ to finding a~limiting distribution for
\begin{multline}\label{betacramerwold}
n^{-1/2}\left([\bX,\bY]^{\top}[\bX,\bY]-\E[\bX,\bY]^{\top}[\bX,\bY]\right)\left[\begin{array}{c}
\bbeta\\-1
\end{array}\right]\\=n^{-1/2}\left([\bI,\bbeta]^{\top}\bZ^{\top}[\btheta,\beps]+[\btheta,\beps]^{\top}[\btheta,\beps]-n\sigma^2\bI\right)\left[\begin{array}{c}
\bbeta\\-1
\end{array}\right].
\end{multline}

Now, it is sufficient to prove the univariate asymptotic normality of
\[
n^{-1/2}\sum_{i=1}^n{\bf t}^{\top}\left([\bZ,\bZ\bbeta]_{i,\bullet}^{\top}[\btheta,\beps]_{i,\bullet}+[\btheta,\beps]_{i,\bullet}^{\top}[\btheta,\beps]_{i,\bullet}-\sigma^2\bI\right)\left[\begin{array}{c}
\bbeta\\-1
\end{array}\right],\quad\forall {\bf t}\in\mathbbm{R}^{p+1};
\]
and apply the Cram\'{e}r-Wold theorem. %If ${\bf t}=c\left[\begin{array}{c}

The case of ${\bf t}=\zero$ is trivial. %On the other hand, if $\zero\neq{\bf t}\neq c\left[\begin{array}{c}
%\bbeta\\-1
%\end{array}\right]$
If ${\bf t}\neq\zero$, %for any real constant $c\neq 0$, 
a~sum of zero mean $\alpha$-mixing random variables
\begin{align}\label{rho_i}
\rho_i&:=\bt^{\top}\left([\bZ,\bZ\bbeta]_{i,\bullet}^{\top}[\btheta,\beps]_{i,\bullet}+[\btheta,\beps]_{i,\bullet}^{\top}[\btheta,\beps]_{i,\bullet}-\sigma^2\bI\right)\left[\begin{array}{c}
\bbeta\\-1
\end{array}\right]\\
&=\bZ_{i,\bullet}\bt_{-(p+1)}\btheta_{i,\bullet}\bbeta+\bZ_{i,\bullet}\bbeta t_{p+1}\btheta_{i,\bullet}\bbeta-\bZ_{i,\bullet}\bt_{-(p+1)}\eps_i-\bZ_{i,\bullet}\bbeta t_{p+1}\eps_i\nonumber\\
&\quad +\btheta_{i,\bullet}\bt_{-(p+1)}\btheta_{i,\bullet}\bbeta+t_{p+1}\eps_i\btheta_{i,\bullet}\bbeta-\btheta_{i,\bullet}\bt_{-(p+1)}\eps_i-t_{p+1}\eps_i^2-\sigma^2\bt^{\top}\left[\begin{array}{c}
\bbeta\\-1
\end{array}\right]\nonumber
\end{align}
satisfies all the assumptions of Corollary~\ref{col:CLTwALPHA}. In fact, assumption~\eqref{WIPwALPHAcond1} holds for $\omega=1/2\min_{j=1,\ldots,p+1}\omega_j$ realizing~\eqref{assalphaZAN}, \eqref{assalpha3AN1}, and~\eqref{jensen-moment} together with the Cauchy-Schwarz inequality and Jensen's inequality:
\begin{align*}
&\sup_{n\in\mathbbm{N}}\E|\rho_n|^{2+\omega}\leq 9^{1+\omega}\sup_{n\in\mathbbm{N}}\E\bigg\{|\bZ_{n,\bullet}\bt_{-(p+1)}\btheta_{n,\bullet}\bbeta|^{2+\omega}+|\bZ_{n,\bullet}\bbeta t_{p+1}\btheta_{n,\bullet}\bbeta|^{2+\omega}\bigg.\displaybreak[0]\\
&\quad \bigg.+|\bZ_{n,\bullet}\bt_{-(p+1)}\eps_n|^{2+\omega}+|\bZ_{n,\bullet}\bbeta t_{p+1}\eps_n|^{2+\omega}+|\btheta_{n,\bullet}\bt_{-(p+1)}\btheta_{n,\bullet}\bbeta|^{2+\omega}\bigg.\displaybreak[0]\\
&\quad\bigg. +|t_{p+1}\eps_n\btheta_{n,\bullet}\bbeta|^{2+\omega}+|\btheta_{n,\bullet}\bt_{-(p+1)}\eps_n|^{2+\omega}+|t_{p+1}\eps_n^2|^{2+\omega} +|\sigma^2(\bt_{-(p+1)}^{\top}\bbeta-t_{p+1})|^{2+\omega}\bigg\}\displaybreak[0]\\
&\leq 9^{1+\omega}\bigg\{\sup_{n\in\mathbbm{N}}|\bZ_{n,\bullet}\bt_{-(p+1)}|^{2+\omega}\sup_{n\in\mathbbm{N}}\E|\btheta_{n,\bullet}\bbeta|^{2+\omega}+\sup_{n\in\mathbbm{N}}|\bZ_{n,\bullet}\bbeta t_{p+1}|^{2+\omega} \sup_{n\in\mathbbm{N}}\E|\btheta_{n,\bullet}\bbeta|^{2+\omega}\bigg.\displaybreak[0]\\
&\quad +\bigg. \sup_{n\in\mathbbm{N}}|\bZ_{n,\bullet}\bt_{-(p+1)}|^{2+\omega} \sup_{n\in\mathbbm{N}}\E|\eps_n|^{2+\omega}  +\sup_{n\in\mathbbm{N}}|\bZ_{n,\bullet}\bbeta t_{p+1}|^{2+\omega} \sup_{n\in\mathbbm{N}}\E|\eps_n|^{2+\omega}\bigg.\displaybreak[0]\\
&\quad \bigg. +\left[\sup_{n\in\mathbbm{N}}\E|\btheta_{n,\bullet}\bt_{-(p+1)}|^{4+2\omega}\right]^{1/2} \left[\sup_{n\in\mathbbm{N}}\E|\btheta_{n,\bullet}\bbeta|^{4+2\omega}\right]^{1/2} \bigg.\displaybreak[0]\\
&\quad \bigg. +\left[\sup_{n\in\mathbbm{N}}\E|t_{p+1}\eps_n|^{4+2\omega}\right]^{1/2} \left[\sup_{n\in\mathbbm{N}}\E|\btheta_{n,\bullet}\bbeta|^{4+2\omega}\right]^{1/2}\bigg.\displaybreak[0]\\
&\quad \bigg. +\left[\sup_{n\in\mathbbm{N}}\E|\btheta_{n,\bullet}\bt_{-(p+1)}|^{4+2\omega}\right]^{1/2} \left[\sup_{n\in\mathbbm{N}}\E|\eps_n|^{4+2\omega}\right]^{1/2} \bigg.\displaybreak[0]\\
&\quad \bigg. +|t_{p+1}|\sup_{n\in\mathbbm{N}}\E|\eps_n|^{4+2\omega}+|\sigma^2(\bt_{-(p+1)}^{\top}\bbeta-t_{p+1})|^{2+\omega} \bigg\}\displaybreak[0]\\
&\leq 9^{1+\omega}\bigg\{p^{2+2\omega}\max_{j=1,\ldots,p}|t_j|^{2+\omega}\sum_{j=1}^p\sup_{n\in\mathbbm{N}}|Z_{n,j}|^{2+\omega}\max_{j=1,\ldots,p}|\beta_j|^{2+\omega}\sum_{j=1}^p\sup_{n\in\mathbbm{N}}\E|\Theta_{n,j}|^{2+\omega} \bigg.\displaybreak[0]\\
&\quad\bigg. +p^{2+2\omega}|t_{p+1}|\max_{j=1,\ldots,p}|\beta_j|^{2+\omega}\sum_{j=1}^p\sup_{n\in\mathbbm{N}}|Z_{n,j}|^{2+\omega}\max_{j=1,\ldots,p}|\beta_j|^{2+\omega}\sum_{j=1}^p\sup_{n\in\mathbbm{N}}\E|\Theta_{n,j}|^{2+\omega} \bigg.\displaybreak[0]\\
&\quad\bigg. +p^{1+\omega}\max_{j=1,\ldots,p}|t_j|^{2+\omega}\sum_{j=1}^p\sup_{n\in\mathbbm{N}}|Z_{n,j}|^{2+\omega}\sup_{n\in\mathbbm{N}}\E|\eps_n|^{2+\omega} \bigg.\displaybreak[0]\\
&\quad\bigg. +p^{1+\omega}|t_{p+1}|\max_{j=1,\ldots,p}|\beta_j|^{2+\omega}\sum_{j=1}^p\sup_{n\in\mathbbm{N}}|Z_{n,j}|^{2+\omega}\sup_{n\in\mathbbm{N}}\E|\eps_n|^{2+\omega} \bigg.\displaybreak[0]\\
&\quad\bigg. +p^{3+2\omega}\max_{j=1,\ldots,p}|t_j|^{2+\omega}\max_{j=1,\ldots,p}|\beta_j|^{2+\omega}\sum_{j=1}^p\sup_{n\in\mathbbm{N}}\E|\Theta_{n,j}|^{4+2\omega} \bigg.\displaybreak[0]\\
&\quad\bigg. +p^{3/2+\omega}|t_{p+1}|^{2+\omega}\max_{j=1,\ldots,p}|\beta_j|^{2+\omega}\left[\sum_{j=1}^p\sup_{n\in\mathbbm{N}}\E|\Theta_{n,j}|^{4+2\omega}\right]^{1/2}\left[\sup_{n\in\mathbbm{N}}\E|\eps_n|^{4+2\omega}\right]^{1/2} \bigg.\displaybreak[0]\\
&\quad\bigg. +p^{3/2+\omega}\max_{j=1,\ldots,p}|t_j|^{2+\omega}\left[\sum_{j=1}^p\sup_{n\in\mathbbm{N}}\E|\Theta_{n,j}|^{4+2\omega}\right]^{1/2}\left[\sup_{n\in\mathbbm{N}}\E|\eps_n|^{4+2\omega}\right]^{1/2} \bigg.\displaybreak[0]\\
&\quad\bigg. +|t_{p+1}|\sup_{n\in\mathbbm{N}}\E|\eps_n|^{4+2\omega}+|\sigma^2(\bt_{-(p+1)}^{\top}\bbeta-t_{p+1})|^{2+\omega} \bigg\}<\infty.
\end{align*}
%because assumption~\eqref{assalphaZAN} and $\sup_{n\in\mathbbm{N}} |Z_{n,j}|^{2+\omega}<\infty$, $\omega>0,\,j\in\{1,\ldots,p+1\}$ are equivalent.

\citet[Theorem~5.2(a)]{Bradley2005} provides $\alpha(\rho_{\circ},n)\leq \alpha(\eps_{\circ},n)+\sum_{j=1}^p\alpha(\btheta_{\circ,j},n)$. Consequently, assumption~\eqref{WIPwALPHAcond2} holds due to the concavity of function $u\mapsto u^{\omega/(2+\omega)}$, $\omega>0$:
\[
\sum_{n=1}^{\infty}\alpha(\rho_{\circ},n)^{\omega/(2+\omega)}\leq\sum_{n=1}^{\infty}\alpha(\eps_{\circ},n)^{\omega/(2+\omega)}+\sum_{j=1}^p\sum_{n=1}^{\infty}\alpha(\btheta_{\circ,j},n)^{\omega/(2+\omega)}<\infty,\quad\omega>0;
\]
%and due to~\eqref{Onalpha1} and~\eqref{Onalpha2}.
which is true because of~\eqref{assalpha4AN} and the fact that
\begin{align*}
\alpha(\btheta_{\circ,j},n)^{\omega/(2+\omega)}&=\mathcal{O}\left(n^{-1-\frac{\delta_j\omega-2}{2+\omega}}\right),\quad \delta_j>2/\omega>0,\,j\in\{1,\ldots,p\};\\
\alpha(\eps_{\circ},n)^{\omega/(2+\omega)}&=\mathcal{O}\left(n^{-1-\frac{\delta_j\omega-2}{2+\omega}}\right),\quad \delta_{p+1}>2/\omega>0.
\end{align*}

Using~\eqref{long-run} and~\eqref{betacramerwold}, let us calculate
\begin{align}\label{varsigma2On}
&\frac{1}{n}\E\left(\sum_{i=1}^n\rho_i\right)^2=\frac{1}{n}\E\left\{\bt^{\top}\left([\bZ,\bZ\bbeta]^{\top}[\btheta,\beps]+[\btheta,\beps]^{\top}[\btheta,\beps]-n\sigma^2\bI\right)\left[\begin{array}{c}
\bbeta\\-1
\end{array}\right]\right\}^2\nonumber\\
&=\frac{1}{n}\E\bt^{\top}\left([\bX,\bY]^{\top}[\bX,\bY]-\E[\bX,\bY]^{\top}[\bX,\bY]\right)\left[\begin{array}{c}
\bbeta\\-1
\end{array}\right]\nonumber\\
&\quad [\bbeta^{\top},-1]\left([\bX,\bY]^{\top}[\bX,\bY]-\E[\bX,\bY]^{\top}[\bX,\bY]\right)\bt\nonumber\\
&=\bt^{\top}\left\{\frac{1}{n}\Var[\bX,\bY]^{\top}[\bX,\bY]\left[\begin{array}{c}
\bbeta\\-1
\end{array}\right]\right\}\bt\nonumber\to\bt^{\top}\beth\bt>0,\quad n\to\infty;
\end{align}
and assumption~\eqref{WIPwALPHAcond3} is satisfied as well. Thus, Corollary~\ref{col:CLTwALPHA} implies asymptotically zero mean normal distribution $n^{-1/2}\sum_{i=1}^n\rho_i$.
\end{Proof}

\begin{Proof}[\emph{Theorem~\ref{thm:ANphi}}]

This proof contains very similar ideas as the proof of Theorem~\ref{thm:ANalpha}.

Assumptions of Theorem~\ref{thm:ANphi} imply the assumptions of Theorem~\ref{thm:CONphi}. Indeed, assumption~\eqref{assphi1a} is implied by assumption~\eqref{assphi1ANa}. Assumption~\eqref{assphi2a} directly follows from~\eqref{assphi2AN1}, because~\eqref{jensen-moment} and \eqref{assphi2AN1} yield
\[
\sup_{n\in\mathbbm{N}}\E\xi_n^4\leq\left(\sup_{n\in\mathbbm{N}}\E|\xi_n|^{4+\omega}\right)^{4/(4+\omega)}<\infty
\]
for $\xi_n\in\{\Theta_{n,1},\ldots,\Theta_{n,p},\eps_n\}$ and $\omega=\min_{j=1,\ldots,p+1}\omega_j$. Hence,
\[
\sum_{n=1}^{\infty}\frac{\E\xi_n^4}{n^2}\leq\sum_{n=1}^{\infty}\frac{\sup_{\iota\in\mathbbm{N}}\E\xi_{\iota}^4}{n^2}=\sup_{\iota\in\mathbbm{N}}\E\xi_{\iota}^4\sum_{n=1}^{\infty}\frac{1}{n^2}<\infty.
\]
Therefore, the consistency results for $\varphi$-mixing errors can be used.

With respect to the Slutsky's theorem, to the Cram\'{e}r-Wold theorem, and to the proof of Theorem~\ref{thm:ANalpha}, it is necessary to find the limiting distribution of $\{\rho_n\}_{n=1}^{\infty}$ from~\eqref{rho_i}. In light of Corollary~\ref{col:CLTwPHI}, we only need to check whether sequences $\{\rho_n\}_{n=1}^{\infty}$ is $\varphi$-mixing sequences, i.e., $\varphi(\rho_{\circ},n)\to 0$ as $n\to\infty$. This follows directly from \citet[Theorem~5.2(d)]{Bradley2005} and assumptions $\varphi(\Theta_{\circ,j},n)\to 0$ for $j=1,\ldots,p$ and $\varphi(\eps_{\circ},n)\to 0$ as $n\to\infty$. The rest of the assumptions of Corollary~\ref{col:CLTwPHI} is included in the assumptions of Corollary~\ref{col:CLTwALPHA} and has been completely checked on sequence $\{\rho_n\}_{n=1}^{\infty}$ in the proof of Theorem~\ref{thm:ANalpha}.
\end{Proof}

%\begin{thebibliography}{000}
% References should be listed in alphabetical order.
\bibliography{tls-weak}

% For books, research reports and proceedings
%\bibitem{B}
%Author(s):
%\newblock{Title of book.}
%\newblock{Edition (if other than first). Publisher's name, place and year of publication (For books, research reports and proceedings).}

% Example
%\bibitem{HoTu96}
%R.~Horst and H.~Tuy:
%\newblock{Global Optimization.}
%\newblock{Springer--Verlag, Berlin 1996.}

% For journal article
%\bibitem{P}
%Author(s):
%\newblock{Title of paper.}
%\newblock{Title of the journal (abbreviated in accordance with Math. Reviews), volume, year of %publication in brackets, inclusive pagination (For journal article).}

% Example
%\bibitem{No85}
%M.~Nov\'{a}k:
%\newblock{A note on the algorithms for determining the model structure.}
%\newblock{Kybernetika {\mi 21} (1985), 164--178.}

% For a paper in a bound collection
%\bibitem{PC}
%Author(s):
%\newblock{Title of paper.}
%\newblock{In: Title of collection, name(s) of editor (in brackets), publisher's name, place and year of publication, inclusive pagination (For a paper in a bound collection).}

% Example
%\bibitem{Pe81}
%V.~Peterka:
%\newblock{Bayesian system identification.}
%\newblock{In: Trends and Progress in System Identification (P.~Eykhoff, ed.), Pergamon Press, Oxford 1981, pp. 239--304.}

%\end{thebibliography}

\makecontacts

\end{document}